\documentclass[12pt,psamsfonts]{amsart}

\newtheorem{anyprop}{Anyprop}

\newtheorem{proposition}[anyprop]{Proposition}

\theoremstyle{definition}

\newtheorem{problem}[anyprop]{Problem}

\newtheorem{remark}[anyprop]{Remark}

\newcommand{\NN}{\mathbb{N}}
\newcommand{\ZZ}{\mathbb{Z}}
\newcommand{\QQ}{\mathbb{Q}}

\newcommand{\CC}{\mathbb{C}}
\newcommand{\FF}{\mathbb{F}}

\newcommand{\PP}{\mathbb{P}}

\newcommand  {\shS}     {\mathcal{S}}
\newcommand  {\shT}     {\mathcal{T}}

\newcommand  {\foB}     {\mathfrak{B}}

\newcommand  {\fom}     {\mathfrak{m}}

\newcommand  {\foo}     {\mathfrak{o}}

\newcommand  {\foX}     {\mathfrak{X}}

\newcommand  {\GL}      {\operatorname{GL}}

\newcommand  {\lra}     {\longrightarrow}

\renewcommand{\O}       {\mathcal{O}}

\newcommand  {\Proj}    {\operatorname{Proj}}

\newcommand  {\ra}      {\rightarrow}

\newcommand  {\rk}    {\operatorname{rk}}

\newcommand  {\Syz}     {\operatorname{Syz}}

\theoremstyle{remark}

\numberwithin{equation}{section}

\usepackage{amscd}
\usepackage{amssymb}

\newcommand{\dwY}{\mathcal Y}

\def\mydate{\number\day\space\ifcase\month \or January\or February\or March\or April\or May\or
June\or July\or August\or September\or October\or November\or
December\fi \space\number\year}

\begin{document}

\title[No restriction for strong
semistability]
{There is no Bogomolov type restriction theorem for strong
semistability in positive characteristic}

\author[Holger Brenner]{Holger Brenner}
\address{Mathematische Fakult\"at, Ruhr-Universit\"at Bochum, 
               44780 Bochum, Germany}
\email{Holger.Brenner@ruhr-uni-bochum.de}





\begin{abstract}
We give an example of a strongly semistable vector bundle of rank two on the projective plane
such that there exist smooth curves of arbitrary high degree with the property
that the restriction
of the bundle to the curve is not strongly semistable anymore.
This shows that a Bogomolov type restriction theorem does not hold
for strong semistability in positive characteristic.
\end{abstract}

\maketitle


\section*{Introduction}

A locally free sheaf $\shS$ on a smooth projective variety $X$ is called semistable if
for every coherent torsion-free subsheaf $\shT \subseteq \shS$ the
inequality $\mu(\shT) \leq \mu(\shS)$ holds, where the slope is defined by
$\mu(\shT )= \deg(\shT) /\rk(\shT)$ and the degree is defined with respect to a fixed ample divisor $H$ on $X$ by $\deg(\shT)= \deg (\det (\shT)) = \det(\shT). H^{\dim X -1}$.
On a curve, this notion is independent of an ample divisor and on a projective
space we always take $H= \O(1)$.

A natural question is whether the restriction of a semistable vector bundle
to a hypersurface is again semistable. One motivation for this is the attempt to reduce
boundedness questions to lower dimensional varieties.
Let us recall briefly what is known
in the easist case,
for bundles on the projective plane and their restrictions to curves
(see \cite[Chapter 7]{huybrechtslehn}).
One kind of result is that the restriction of a semistable bundle to a generic curve of sufficiently high degree is again semistable.
The theorem of Flenner (\cite{flennerrestriction}) gives an explicit
bound for the degree of the curve in characteristic zero;
the theorem of Mehta-Ramanathan (\cite{mehtaramanathanrestriction}) holds in arbitrary characteristic, but without a degree bound.

The theorem of Bogomolov (\cite{bogomolovstability}) shows in characteristic zero that the restriction
of a stable bundle is again  stable for every smooth curve of sufficiently high degree,
not only for the generic curve.
A. Langer has recently shown that such a restriction theorem holds also in positive characteristic, see \cite[Theorem 5.2]{langersemistable}.

A well known phenomenon in positive characteristic is that the pull-back of a semistable
sheaf under the Frobenius morphism is in general no longer semistable.
A bundle is called strongly semistable if every Frobenius pull-back of it
is again semistable. A semistable bundle on the projective plane is also strongly semistable
(the same is true for elliptic curves).
What is known about restriction of strong semistability?
A. Langer has shown in positive characteristic that under the condition that
the discriminant of the bundle is $0$, then its restriction to a very general curve is strongly semistable \cite[Theorem 3.1 and Remark 3.10.2]{langersemistable}.

The general question how strong semistability behaves under restriction is open.
In this paper we show by an example that we cannot hope for the best possible result,
that is a result \`{a} la Bogomolov does not hold: there is no restriction theorem
for strong semistability which holds for all smooth curves of sufficiently high degree.

Our example is motivated by the theory of tight closure and quite easy to write down, it is the syzygy bundle $\Syz(X^a,Y^a,Z^a)$ on the projective plane
(which is the cotangent bundle for $a=1$). We show that for every $d_0$ there exists a
$d \geq d_0$
(depending on the characteristic $p$) such that the restriction to the smooth
Fermat curve $X^d+Y^d+Z^d=0$ is not strongly semistable.

The content of this paper is as follows:
We describe the relevant properties of the example in section \ref{example}.
In section \ref{remarks} we make several remarks concerning further implications of the example,
in particular with respect to tight closure.
We show that $(YXZ)^b \in (X^{2b},Y^{2b},Z^{2b})^*$ (which is always true in zero
characteristic for solid closure) does not hold
in the Fermat hypersurface domains $R=K[X,Y,Z]/(X^d+Y^d+Z^d)$
for infinitely many degrees $d$ (Remark \ref{tight}).
In section \ref{werner}
we apply the example to the $p$-adic setting of a preprint of C. Deninger and A. Werner
(\cite{deningerwerner}) and show that their category
$\foB_{\foX_{\foo}} $ of vector bundles cannot contain
all irreducible bundles of degree $0$.

I thank A. Langer and A. Werner for useful discussion.

\section{The example}
\label{example}

We describe now the example.
Let $p$ denote a prime number and let $K$ denote a field of characteristic $p$.
We consider the syzygy bundle $\Syz(X^a,Y^a,Z^a)(n)$, $a \geq 1$, on $\PP^2= \Proj K[X,Y,Z]$
given by the defining sequence
$$0 \lra \Syz(X^a,Y^a,Z^a)(n)     \lra \bigoplus_3 \O(n-a)
\stackrel{X^a,Y^a,Z^a}{\lra } \O(n) \lra 0 \, .$$
For $a=1$ and $n=0$ this is just the Euler-Sequence for the cotangent bundle.
The pull-back of this sequence under the $e$-th absolute Frobenius morphism
$F^{e}:\PP^2 \lra \PP^2$ yields the sequences (setting $q=p^{e}$)
$$0 \lra (F^{e*}\Syz(X^a,Y^a,Z^a))(m)  \lra \bigoplus_3 \O(m-aq)
\stackrel{X^{aq},Y^{aq},Z^{aq}}{\lra } \O(m) \lra 0 \, .$$
Therefore we have $(F^{e*} \Syz(X^a,Y^a,Z^a)) (m) = \Syz(X^{aq},Y^{aq},Z^{aq}) (m)$.
For this bundle the following proposition holds.

\begin{proposition}
\label{strongprop}
Let $K$ denote a field of positive characteristic $p$.
Then the syzygy bundle $\Syz(X^a,Y^a,Z^a)$ on the projective plane $\PP^2$
is strongly semistable but it is not true
that the restriction of this bundle is strongly semistable
for all smooth curves of sufficiently high degree.
\end{proposition}

\proof
First note that the bundles $\Syz(X^{aq},Y^{aq},Z^{aq})(m)$
do not have non-trivial global sections for $m <2aq$ for all $a$ and $q$.
So these bundles are semistable and also strongly semistable on the projective plane.
For the main statement we have to show that for every number $d_0 \in \NN$ we find
a degree $d \geq d_0$ and a smooth projective curve $C \subseteq  \PP^2$
of degree $d$ such that the restriction
$\Syz(X^a,Y^a,Z^a)|C$ is not strongly semistable anymore.

Let $d_0$ be given and let $e$ be such that
$ ap^{e-1} \geq d_0$.
Choose $d$ such that
$$ap^{e-1} < d < 3ap^{e-1}/2$$
and that $d$ is no multiple of $p$.
Then $$aq=ap^{e} < dp < 3aq/2 \, .$$
Set $k= dp-aq >0$.
Consider the smooth curve $C$ given by the Fermat equation
$X^d +Y^d +Z^d=0$.
On the affine cone $K[X,Y,Z]/(X^d +Y^d +Z^d)$ over this curve we have
$$0= (X^d+Y^d+Z^d)^p = X^{dp} +Y^{dp} +Z^{dp} =
X^k X^{aq}+Y^kY^{aq} +Z^kZ^{aq} \, .$$
This means that on the curve $C$ the syzygy bundle
$\Syz (X^{aq},Y^{aq},Z^{aq})(m)$ has a global non-trivial
section given by $(X^k,Y^k,Z^k)$
of total degree $m= k+aq= dp$.
But the degree of this syzygy bundle is by the short exact sequence
$$\deg( \Syz(X^{aq},Y^{aq},Z^{aq})(m))=
2m-3aq = 2dp-3aq <0$$
because of the second condition in the choice of $d$.
So this bundle on the curve has negative degree but it has a global section
$\neq 0$. Hence it is not semistable and $\Syz(X^a,Y^a,Z^a)|C$
is not strongly semistable.
\qed

\section{Some remarks}
\label{remarks}

\begin{remark}
\label{HN}
Let $aq=ap^e < dp <3aq/2$ and $k=dp-aq$ as in the proof of
Proposition \ref{strongprop}.
The global section in the restricted syzygy bundle
gives the following
short exact sequence on the Fermat curve $C=V_+(X^d+Y^d+Z^d) \subset \PP^2$: 
$$0 \lra \O_C \lra \Syz(X^{aq},Y^{aq},Z^{aq})(dp) \lra \O_C(2dp-3aq) \lra 0 \, .$$
This is of course the Harder-Narasimhan filtration of $\Syz(X^{aq},Y^{aq},Z^{aq})(dp)$,
and the maximal slope of $F^{e*} \Syz(X^a,Y^a,Z^a)$ is $0$ and its minimal slope is
$(2dp-3aq)d $, since $d$ is the degree of $\O_C(1)$.

A measure for the deviation from strong semistability of a locally free sheaf $\shS$
is given by the number
$$\lim_{e \in \NN} \frac{ \mu_{\rm max} (F^{e*} \shS ) - \mu_{\rm min} (F^{e *} \shS)}{p^{e}} \,.$$
This number is independent of a twist with an invertible sheaf. Thus in our example
this number is
$\geq (0 -(2dp-3aq))d/q = (3aq-2dp)d/q$.
If we set $d=ap^{e-1}+1$, then
this number is
$$ \frac{d(3aq -2p(ap^{e-1} +1))}{q} =  \frac{d(aq-2p)}{q}
\geq \frac{ap^{e-1}(aq-2p)}{q} = a^2p^{e-1} -2a  \, .$$
So for every $a$ this deviation measure is arbitrary big (but finite if we divide by $d$).
\end{remark}

\begin{remark}
As mentioned in the introduction,
there exist restriction theorems due to Bogomolov and to
A. Langer for semistable sheaves
which imply that the restriction to every smooth hypersurface of sufficiently high degree is
again semistable. Here the degree bound does not depend only on the variety and the rank of the bundle (as in the restriction theorem of Flenner),
but also on the discriminant of the bundle.
The example above shows that in positive characteristic there cannot exist a restriction theorem
for semistability for all smooth hypersurfaces
with a degree bound depending only on data which are invariant under the Frobenius
(like the rank and invariants of the basic variety).
For otherwise the restriction of a strongly semistable sheaf
to a smooth hypersurface would also be strongly semistable.
\end{remark}

\begin{remark}
\label{tight}
The example is motivated by the geometric interpretation of the theory of tight closure.
We refer to \cite{hunekeapplication}
and \cite{hunekeparameter} for background of this theory.
In \cite[Theorem 8.4]{brennerslope} we proved that the tight closure of the ideal
$(X^a,Y^a,Z^a)$ in $R=K[X,Y,Z]/(F)$, where $F$ is homogeneous,
is given by
$$(X^a,Y^a,Z^a)^* =(X^a,Y^a,Z^a) + R_{\geq 3a/2} \hspace{2cm} (*)$$
under the condition that the syzygy bundle
$\Syz(X^a,Y^a,Z^a)$ is strongly semi\-stable on the curve $C=\Proj \,R$.
Since this syzygy bundle is semi\-stable on the projective plane,
it follows in characteristic zero that this numerical formula for tight closure (or
rather solid closure) holds
in $R=K[X,Y,Z]/(F)$ for every $F$ of sufficiently high degree by the restriction theorem of Bogomolov. In fact this holds for every $F$ smooth of degree $\deg (F)=d \geq 3a-1$,
see \cite[Proposition 6.2]{brennercomputationtight}.

Here we want to remark that this is not true in positive characteristic.
Let $p$ be a prime number. We show that for $a=2b$ even
the formula does not hold for $(X^a,Y^a,Z^a)$ on $F=X^d+Y^d+Z^d$ for suitable but arbitrary
high degree $d$;
in fact we will show that $X^bY^bZ^b \not\in (X^a,Y^a,Z^a)^*$
on such a curve. Even degree has the (dis-)advantage that there exists
homogeneous elements having the critical degree $3a/2$.

Suppose as before that $aq< dp< 3aq/2$ and $k=dp-aq$.
From the Harder-Narasimhan filtration (see Remark \ref{HN})
we get by tensoring with $\O(m-dp)$
the short exact sequence
$$0 \lra \O_C(m-k-aq) \lra \Syz(X^{aq},Y^{aq},Z^{aq})(m) \lra \O_C(m+k-2aq) \lra 0 \, .$$
For $4aq-2k > 2m \geq 3aq$ the invertible sheaf on the right has negative degree
but the degree in the middle is $\geq 0$.
Let $f=X^rY^sZ^t$ be a monomial of degree $r+s+t=m$. Such an element yields via the defining
sequence for $\Syz(X^{aq},Y^{aq},Z^{aq})$ the $\check{\rm C}$ech
cohomology class
$$ (\frac{X^rY^sZ^t}{X^{aq}} , -\frac{X^rY^sZ^t}{Y^{aq}}, 0)
 \in H^1(C,\Syz(X^{aq},Y^{aq},Z^{aq})(m)) \, .$$
The surjection
$\Syz(X^{aq},Y^{aq},Z^{aq})(m) \lra \O_C(m+k-2aq)$
maps this cohomology class to the class
$$-\frac{X^rY^sZ^tZ^k}{X^{aq}Y^{aq}} \in H^1(C, \O_C(m+k-2aq))\, ,$$
see \cite[Remark 4.8]{brennercomputationtight}.

If the formula ($*$) would hold,
then $X^rY^sZ^t \in (X^{aq},Y^{aq},Z^{aq})^*$, since $m \geq 3aq/2$.
If however the cohomology class $\frac{X^rY^sZ^tZ^k}{X^{aq}Y^{aq}}$ is not in $0^*$,
which means by definition that $X^rY^sZ^tZ^k \not\in (X^{aq},Y^{aq})^*$ in
$R=K[X,Y,Z]/(F)$,
then (the corresponding geometric torsor of this class is affine and this is then true for
the torsor corresponding to the first class) also
$X^rY^sZ^t \not\in (X^{aq},Y^{aq},Z^{aq})^*$.

Let now be $a=2b \geq 2$, $d$ and $k$ as before and set $m=3aq/2 = 3bq$.
Set $r=s=t=bq$, so that we have to look at the cohomology class
$$\frac{X^{bq}Y^{bq}Z^{bq}Z^k}{X^{aq}Y^{aq}} =\frac{Z^{bq}{Z^k}}{X^{bq}Y^{bq}} \in H^1(C,\O_C(k-bq)) \, ,$$
and we have to look for examples such that
$Z^{bq+k} \not\in (X^{bq}, Y^{bq})^*$.

We fix now the degree $d=ap^{e-1}+1$, so that $k=p$ and $d >p$. Set $u= \lceil p/2 \rceil$.
Then for $p \geq 3$ and $e \geq 2$ we have
$$ ud= \lceil p/2 \rceil (2bp^{e-1} +1)=    bp^{e} + p/2 +1/2 + bp^{e-1} \geq bq+p  \, .$$
So it is enough to show that $Z^{ud} \not\in (X^{bq},Y^{bq})^*$.
But this class is
$$-\frac{Z^{ud}}{X^{bq}Y^{bq}}= \frac{  (X^d+Y^d)^u}{X^{bq}Y^{bq}}
= \frac{ X^{ud}+uX^{(u-1)d}Y^d + \ldots + Y^{ud} }{X^{bq}Y^{bq}}\, ,$$
and we consider this class in $H^1(\PP^1, \O(ud-2bq))$, $\PP^1=\Proj\, K[X,Y]$.
The extreme summands on the left and on the right vanish,
but for $e \geq 2$ and $p \geq 5$ we got that
$(u-1)d < bq$, so the second summand is $\neq 0$ ($p$ does not divide $u$).
Therefore this class is $\neq 0$ and
hence $\not\in 0^*$, since $K[X,Y]$ is $F$-regular.
It follows that the cohomology class does not belong to
$0^* \subseteq H^1(C, \O_C(k-bq))$
and so $(XYZ)^b  \not\in (X^{2b},Y^{2b},Z^{2b})^*$ in $R=K[X,Y,Z]/(F)$.
\end{remark}

\begin{problem}
Describe in dependence of the prime number $p$ and the degree $d$
whether $\Syz(X^2,Y^2,Z^2)$
is strongly semistable on $X^d+Y^d+Z^d=0$.
In particular, for fixed degree $d$,
decide whether the set of prime numbers such that this syzygy bundle is strongly semistable
is finite, infinite or contains almost all prime numbers?
This question is both a special instance of a problem of Miyaoka
(see \cite[Problem 5.4]{miyaokachern}) and of a problem of tight closure theory
(\cite[Question 13]{hochstertightsolid}).
\end{problem}

\section{How big are $\foB_{X_{\CC_p}}$ and $\foB_{\foX_{\foo}}$?}
\label{werner}

C. Deninger and A. Werner establish in \cite{deningerwerner}
a $p$-adic analog of the well-known analytic correspondence
between flat vector bundles over a compact Riemann surface $X$
and $\GL$-representations of the geometric fundamental group $\pi_1(X,x)$.
They construct a functor which associates to a certain category of vector
bundles on a smooth projective $p$-adic curve $X_{\CC_p}$ a continuous
representation of the algebraic fundamental group in a vector space over $\CC_p$.
However it remains open how big this category of bundles is,
in particular it is unclear whether it contains all irreducible (or all stable) bundles
of degree $0$.

We describe briefly the setting of \cite{deningerwerner}.
Let $K$ denote a finite extension of the $p$-adic number field $\QQ_p$,
and let $ \overline{K}$ denote its algebraic closure.
Let $\CC_p$ denote the completion of $\overline{K}$
with respect to the p-adic valuation.
Let $X$ denote a smooth projective curve over $K$
and set $\overline{X} =X \otimes _K \overline{K}$ and $X_{\CC_p}=X \otimes_K \CC_p$.
Let $\overline{\foX}$ denote a finitely presented flat and proper model over
the valuation domain $\foo_{\overline{K}} \subset \overline{K}$ so that $\overline{X}$ is the generic fiber.
The special fiber is then a not necessarily irreducible
projective curve over $\foo_{\overline{K}}/\fom = \overline{\FF}_p$

They work with a special category
$\foB= \foB_{\foX_\foo}$ of vector bundles on $\foX_\foo = \overline{\foX} \otimes \foo$,
where $\foo$
denotes the valuation ring in $\CC_p$.
A vector bundle $E$ on $\foX_{\foo}$
belongs to this category if for every $n \geq 1$
there is a (generic \'{e}tale and proper) ``covering''
$\pi: \dwY \ra \overline{\foX}$ (with some further technical properties)
such that $\pi^* E_n$ ($ E_n =E/(p^n)$) is trivial on the thickened fiber
$\dwY_n =\dwY/(p^n)$.

For the bundles in this category
they construct a functor
$\rho$ from
$\foB_{\foX_\foo}$ to the continuous representations
${\bf Rep}_{\pi_1( \overline{X}, \overline{x})} (\foo)$
of the algebraic fundamental group $\pi_1(\overline{X}, \overline{x})$
(\cite[Theorem 14]{deningerwerner}).
This functor yields then also a functor (\cite[Theorem 23]{deningerwerner})
$$\rho: \foB_{X_{\CC_p}} \lra  {\bf Rep}_{\pi_1(\overline{X},\overline{x})} (\CC_p) \, ,$$
where the category  $\foB_{X_{\CC_p}}$
consists of the vector bundles
on $X_{\CC_p}$ which arise by restriction from a bundle in $\foB_{\foX_\foo}$
for some model $\overline{\foX}$ of $\overline{X}$ (\cite[Definition 19]{deningerwerner}).
The category $\foB_{X_{\CC_p}}$ contains all invertible sheaves of degree $0$
(\cite[Theorem 23c]{deningerwerner})
and it is hoped that it contains all irreducible bundles of degree
$0$ (\cite[Question 7d]{deningerwerner}). Our example shows that this is definitely not true for $\foB_{\foX_{\foo}}$!

Since the Fermat equation
$U^d+Y^d+Z^d=0$ is defined over $\ZZ$ we may consider
the corresponding projective curve over any ring.
We denote this curve over $\overline{\QQ}_p =\overline{K}$ by $\overline{X}$ and the relative curve over $\foo_{\overline{K}}$ 
by $\overline{\foX}$, so this is our $p$-adic model for $\overline{X}$.
We consider the bundle
$E= \Syz(U^2,Y^2,Z^2)(3)$ on $\foX_{\foo}$.
The determinant of this bundle is just the structure sheaf.
This bundle has degree $0$ and is stable and irreducible on $\overline{X}$,
on $X_{\CC_p}$ and on the special fiber $\foX_\foo \otimes \foo/\fom$,
since this holds over every field (of any characteristic) for $d$ big enough.

The condition in the definition of $\foB_{\foX_\foo}$,
namely that there exists $\pi: \dwY \ra \overline{\foX}$
such that $\pi^* (E_n)$ is trivial,
implies in particular that the bundle over the special fiber
triviales in a finite extension.
A result of Lange and Stuhler ([Satz 1.9]\cite{langestuhler})
however says that a vector bundle on a smooth projective
curve over an algebraically closed field of positive characteristic
must be strongly semistable if it admits such a finite trivialization.
But due to Proposition \ref{strongprop} this syzygy bundle is not strongly semistable for
some degree $d$.
Therefore for suitable degree $d$ the bundle
$E= \Syz(U^2,Y^2,Z^2)(3)$
does not belong to $\foB_{\foX_{\foo}}$.

It is however possible that there exists
another model $\foX'$ and a vector bundle
$E' \in \foB_{\foX'_{\foo}}$
such that $E'_{ \CC_p} \cong E_{\CC_p}$.
Then $E_{\CC_p} \in \foB_{X_{\CC_p}}$ would hold.

\bibliographystyle{plain}

\bibliography{bibliothek}

\end{document}